\documentclass[a4paper,11pt]{article}
\usepackage{amsfonts}
\usepackage{amstext}
\usepackage{amsthm}
\usepackage{amsmath}
\usepackage{amssymb}
\usepackage{graphicx}
\usepackage{latexsym}
\usepackage[latin1]{inputenc}
\newcommand{\R}{\Bbb{R}}
\newcommand{\N}{\Bbb{N}}

\newtheorem{teor}{Theorem}[section]
\newtheorem{propo}{Proposition}[section]

\newcommand{\n}{\noindent}

\newcommand {\fim}{\rule{0.5em}{0.5em}}

\begin{document}

\title{On strongly indefinite systems involving
\\fractional elliptic operators
}

\author{\textbf{Edir Junior Ferreira Leite \footnote{\textit{E-mail addresses}:
edirjrleite@ufv.br (E.J.F. Leite)}}\\ {\small\it Departamento de Matem\'{a}tica,
Universidade Federal de Vi\c{c}osa,}\\ {\small\it CCE, 36570-000, Vi\c{c}osa, MG, Brazil}
}
\date{}{

\maketitle

\markboth{abstract}{abstract}
\addcontentsline{toc}{chapter}{abstract}

\hrule \vspace{0,2cm}

\n {\bf Abstract}

In this paper we discuss the existence and regularity of solutions of the following strongly indefinite systems involving fractional elliptic operators on a smooth bounded domain $\Omega$ in $\R^n$:
\[
\left\{
\begin{array}{llll}
\mathcal{L}u = \mu v+v^p & {\rm in} \ \ \Omega\\
\mathcal{L}v = \lambda u+u^q & {\rm in} \ \ \Omega\\
u= v=0 & {\rm on} \ \ \Sigma
\end{array},
\right.
\]
where $\mathcal{L}$ refer to any of the two types of operators $\mathcal{A}^s$ or $(-\Delta)^s$, $0 < s < 1$, $p,q>1$, $\lambda$ and $\mu$ are fixed real numbers and

\begin{itemize}
      \item[$\bullet$] $\Sigma = \partial\Omega$ for the spectral fractional Laplace operator $\mathcal{A}^s$,
      \item[$\bullet$] $\Sigma = \mathbb{R}^n\setminus\Omega$ for the restricted fractional Laplace operator $(-\Delta)^s$.
\end{itemize}

\vspace{0.5cm}
\hrule\vspace{0.2cm}

\section{Introduction and main result}

This work is devoted to the study of existence of solutions for nonlocal elliptic systems on bounded domains which will be described henceforth.

The fractional Laplace operator (or fractional Laplacian) of order $2s$, with $0 < s < 1$, denoted by $(-\Delta)^{s}$, is defined as

\[
(-\Delta)^{s}u(x) = C(n,s)\, {\rm P.V.}\int\limits_{\R^{n}}\frac{u(x)-u(y)}{\vert x-y\vert^{n+2s}}\; dy\, ,
\]
or equivalently,

\[
(-\Delta)^{s}u(x)=-\frac{1}{2}C(n,s)\int\limits_{\R^{n}}\frac{u(x+y)+u(x-y)-2u(x)}{|y|^{n+2s}}\; dy
\]

\n for all $x \in \R^{n}$, where P.V. denotes the principal value of the integral and

\[
C(n,s) = \left(\int\limits_{\R^{n}}\frac{1-\cos(\zeta_{1})}{\vert\zeta\vert^{n+2s}}\; d\zeta\right)^{-1}
\]

\n with $\zeta = (\zeta_1, \ldots, \zeta_n) \in \R^n$.

We remark that $(-\Delta)^{s}$ is a nonlocal operator on functions compactly supported in $\R^n$. The convergence property

\[
\lim_{s\rightarrow 1^{-}}(-\Delta)^{s} u = -\Delta u
\]

\n pointwise in $\R^n$ holds for every function $u \in C^{\infty}_{0}(\R^{n})$, so that the operator $(-\Delta)^{s}$ interpolates the Laplace operator in $\R^n$.

Factional Laplace operators arise naturally in several different areas such as Probability, Finance, Physics, Chemistry and Ecology, see \cite{A, bucur}.

A closely related operator but different from $(-\Delta)^{s}$, the spectral fractional Laplace operator $\mathcal{A}^{s}$, is defined in terms of the Dirichlet spectra of the Laplace operator on $\Omega$, see \cite{compare, rafaella}. Roughly, if $(\varphi_k)$ denotes a $L^2$-orthonormal basis of eigenfunctions corresponding to eigenvalues $(\lambda_k)$ of the Laplace operator with zero Dirichlet boundary values on $\partial \Omega$, then the operator $\mathcal{A}^s$ is defined as $\mathcal{A}^{s} u = \sum_{k=1}^\infty c_k \lambda_k^s \varphi_k$, where $c_k$, $k \geq 1$, are the coefficients of the expansion $u = \sum_{k=1}^\infty c_k \varphi_k$.

After the work \cite{CaSi} on the characterization for any $0 < s < 1$ of the operator $(-\Delta)^{s}$ in terms of a Dirichlet-to-Neumann map associated to a suitable extension problem, a great deal of attention has been dedicated in the last years to nonlinear nonlocal problems of the kind

\begin{equation}\label{3}
\left\{
\begin{array}{rrll}
(-\Delta)^{s} u &=& f(x,u) & {\rm in} \ \ \Omega\\
u &=& 0  & {\rm in} \ \ \R^n\setminus\Omega
\end{array}
\right.
\end{equation}
where $\Omega$ is a smooth bounded open subset of $\R^{n}$, $n \geq 1$ and $0 < s < 1$.

Several works have focused on the existence \cite{ros131, ros137, ros138, auto2, auto1, ros215, niang, ros264, ros265, ros267, val}, nonexistence \cite{ros129, val}, symmetry \cite{ros21, ros95} and regularity \cite{ros10, cabre1, ROS1, ROS, S} of viscosity solutions, among other qualitative properties \cite{ros1, CR, ros130}.

A specially important example arises for the power function $f(x,u) = u^p$, with $p > 0$, in which case (\ref{3}) is called the fractional Lane-Emden problem. Recently, it has been proved in \cite{ros264} that this problem admits at least one positive viscosity solution for $1 < p < \frac{n + 2s}{n - 2s}$. The nonexistence has been established in \cite{ROS2} whenever $p \geq \frac{n + 2s}{n - 2s}$ and $\Omega$ is star-shaped. These results were known long before for $s = 1$, see the classical references \cite{AR, GS, ros237, Rabinowitz} and the survey \cite{Pucci}.

An extension for spectral fractional operator was devised by Cabr\'{e} and Tan \cite{CT} and Capella, D\'{a}vila, Dupaigne, and Sire \cite{CDDS} (see Br\"{a}ndle, Colorado, de Pablo, and S\'{a}nchez \cite{BrCPS} and Tan \cite{T} also). Thanks to these advances, the boundary fractional problem

\begin{equation} 
\left\{
\begin{array}{rrll}
\mathcal{A}^{s} u &=& u^p & {\rm in} \ \ \Omega\\
u &=& 0  & {\rm on} \ \ \partial\Omega
\end{array}
\right.
\end{equation}

\n has been widely studied on a smooth bounded open subset $\Omega$ of $\mathbb{R}^n$, $n\geq 2$, $s \in (0,1)$ and $p > 0$. Particularly, a priori bounds and existence of positive solutions for subcritical exponents ($p < \frac{n + 2s}{n - 2s}$) has been proved in \cite{BrCPS, CT, choi, CK, T} and nonexistence results has also been proved in \cite{BrCPS, tan1, T} for critical and supercritical exponents ($p \geq \frac{n + 2s}{n - 2s}$). The regularity result has been proved in \cite{cabre1, CDDS, T, yang}.

When $s = 1/2$, Cabr\'{e} and Tan \cite{CT} established the existence of positive solutions for equations having nonlinearities with the subcritical growth, their regularity, the symmetric property, and a priori estimates of the Gidas-Spruck type by employing a blow-up argument along with a Liouville type result for the square root of the Laplace operator in the half-space. Then \cite{T} has the analogue to $1/2<s<1$. Br\"{a}ndle, Colorado, de Pablo, and S\'{a}nchez \cite{BrCPS} dealt with a subcritical concave-convex problem. For $f(u)=u^q$ with the critical and supercritical exponents $q\geq\frac{n+2s}{n-2s}$, the nonexistence of solutions was proved in \cite{BaCPS, tan1, T} in which the authors devised and used the Pohozaev type identities. The Brezis-Nirenberg type problem was studied in \cite{tan1} for $s = 1/2$ and \cite{BaCPS} for $0<s<1$. The Lemma's Hopf and Maximum Principe was studied in \cite{T}.

An interesting interplay between the two operators occur in case of periodic solutions, or when the domain is the torus, where they coincide, see \cite{torus}. However in the case general the two operators produce very different behaviors of solutions, even when one focuses only on stable solutions, see e.g. Subsection 1.7 in \cite{serena5}.

We here are interested in studying the following problem

\begin{equation}\label{1}
\left\{
\begin{array}{llll}
\mathcal{L}u = \mu v+v^p & {\rm in} \ \ \Omega\\
\mathcal{L}v = \lambda u+u^q & {\rm in} \ \ \Omega\\
u= v=0 & {\rm on} \ \ \Sigma
\end{array},
\right.
\end{equation}
where $\mathcal{L}$ refer to any of the two types of operators $\mathcal{A}^s$ or $(-\Delta)^s$, $0 < s < 1$, $p,q>1$, $\lambda$ and $\mu$ are fixed real numbers and

\begin{itemize}
      \item[$\bullet$] $\Sigma = \partial\Omega$ for the spectral fractional Laplace operator $\mathcal{A}^s$,
      \item[$\bullet$] $\Sigma = \mathbb{R}^n\setminus\Omega$ for the restricted fractional Laplace operator $(-\Delta)^s$.
\end{itemize}

For $0<s\leq 1$ and $\lambda=0=\mu$, the problem (\ref{1}) and a number of its generalizations have been widely investigated in the literature during the two last decades. For $s=1$ see for instance \cite{CFM, FM, DG, DF, DR, vander, Mi1, Mi2, marcos, SZ1} for $\lambda=0=\mu$ and see \cite{vander} for $\lambda$ and $\mu$ fixed real numbers. Now for $0<s<1$ and $\lambda=0=\mu$ the system above was investigated in \cite{EM1} for $\mathcal{L}=(-\Delta)^s$ and in \cite{choi} for $\mathcal{L}=\mathcal{A}^s$.

Related systems have been investigated by using other methods. We refer to the works \cite{serena2, serena1, EM} for systems involving different operators $(-\Delta)^{s}$ and $(-\Delta)^{t}$ in each one of equations. More generally, fractional systems have been studied with extension methods in \cite{serena3, serena4}.

In this work we discuss existence and regularity of solutions of problem (\ref{1}) for $0 < s < 1$, $\lambda$ and $\mu$ fixed real numbers and $\mathcal{L}=(-\Delta)^s$ or $\mathcal{L}=\mathcal{A}^s$. We determine the precise set of exponents $p$ and $q$ for which the problem (\ref{1}) admits always a solution.

The ideas involved in our proofs base on variational methods. In particular, we obtain existence result for Hamiltonian systems through Indefinite Functional Theorem of Benci and Habinowitz. 

Our main result is

\begin{teor}\label{existhamiltinian}
Suppose that $\lambda$ and $\mu$ are fixed real numbers and $p,q>1$ satisfies 
\begin{equation}\label{5}
\frac{1}{p + 1} + \frac{1}{q + 1} > \frac{n - 2s}{n}\, .
\end{equation}
Let $0<\alpha<2s$ be such that
\[
\frac{1}{2}-\frac{1}{q+1}<\frac{\alpha}{n}\text{ and }\frac{1}{2}-\frac{1}{p+1}<\frac{2s-\alpha}{n}.
\]
Then, there exists a (nontrivial) weak solution $(u,v)\in \Theta^{\alpha}(\Omega)\times\Theta^{2s-\alpha}(\Omega)$ to the problem (\ref{1}). Moreover, there is $0<\eta<1$ such that $(u,v)\in (C^{\eta}(\mathbb{R}^n))^2$ if $\mathcal{L}=(-\Delta)^s$ and $(u,v)\in (C^{2,\eta}(\Omega)\cap C^{1,\eta}(\overline{\Omega}))^2$ if $\mathcal{L}=\mathcal{A}^s$.
\end{teor}

\begin{figure}[ht]
\centering
\includegraphics[scale=1.11]{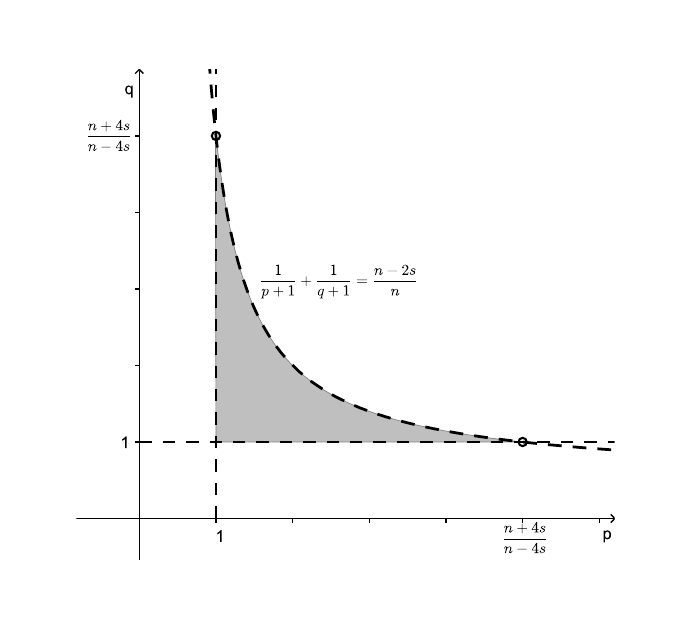}
\caption{The existence range of couples $(p,q)$ when $n>4s$.}
\end{figure}

The rest of paper is organized into two sections. In Section 2 we briefly recall some definitions and facts dealing with fractional Sobolev spaces and comment some relationships and differences between operators $(-\Delta)^{s}$ and $\mathcal{A}^s$. In Section 3, we prove Theorem \ref{existhamiltinian} by applying the Indefinite Functional Theorem of Benci and Rabinowitz. Finally we shall establish the Brezis-Kato type result and study the regularity of solutions to (\ref{1}).

\section{Preliminaries}

We start by fixing a parameter $0 < s < 1$. Let $\Omega$ be an open subset of $\R^n$, with $n \geq 1$. For any $r \in (1,+\infty),$ one defines the fractional Sobolev space $W^{s,r}(\Omega)$ as

\begin{equation}\label{frac2.1}
W^{s,r}(\Omega):=\left\{u\in L^r(\Omega): \frac{\vert u(x)-u(y)\vert}{\vert x-y\vert^{\frac{n}{r}+s}}\in L^r(\Omega\times\Omega)\right\}\, ,
\end{equation}
that is, an intermediary Banach space between $L^r(\Omega)$ and $W^{1,r}(\Omega)$ induced with the norm
\begin{equation}\label{frac2.2}
\Vert u\Vert_{W^{s,r}(\Omega)}:=\left(\int\limits_{\Omega}\vert u\vert^r dx + \int\limits_{\Omega}\int\limits_{\Omega}\frac{\vert u(x)-u(y)\vert^r}{\vert x-y\vert^{n+sr}}dxdy\right)^{\frac{1}{r}}\, ,
\end{equation}
where the term
\[
[u]_{W^{s,r}(\Omega)}:=\left(\int\limits_{\Omega}\int\limits_{\Omega}\frac{\vert u(x)-u(y)\vert^r}{\vert x-y\vert^{n+sr}}dxdy\right)^{\frac{1}{r}}
\]
is the so-called Gagliardo semi-norm of $u$.

Let $s \in \R \setminus \N$ with $s\geq 1$. The space $W^{s,r}(\Omega)$ is defined as
\[
W^{s,r}(\Omega)=\{u\in W^{[s],r}(\Omega) : D^ju\in W^{s-[s],r}(\Omega),\forall j, \vert j\vert=[s]\}\, ,
\]

\n where $[s]$ is the largest integer smaller than $s$, $j$ denotes the $n$-uple $(j_1, \ldots, j_n) \in \N^n$ and $|j|$ denotes the sum $j_1 + \ldots + j_n$.

It is clear that $W^{s,r}(\Omega)$ endowed with the norm

\begin{equation}\label{frac2.11}
\Vert u\Vert_{W^{s,r}(\Omega)}=\left(\Vert u\Vert^r_{W^{[s],r}(\Omega)}+ [u]^r_{W^{s-[s],r}(\Omega)}\right)^{\frac{1}{r}}
\end{equation}
is a reflexive Banach space.

Clearly, if $s=m$ is an integer, the space $W^{s,r}(\Omega)$ coincides with the Sobolev space $W^{m,r}(\Omega)$.

Let $W^{s,r}_0(\Omega)$ denote the closure of $C^\infty_0(\Omega)$ with respect to the norm $\Vert\cdot\Vert_{W^{s,r}(\Omega)}$ defined in (\ref{frac2.11}). For $0 < s \leq 1$, we have
\[
W^{s,r}_0(\Omega)=\{u\in W^{s,r}(\R^n) : u=0 \text{ in }\R^n\setminus\Omega\}\, .
\]
For more details on the above claims, we refer to \cite{DD, frac, tartar}.\\

In this paper, we focus on the case $r = 2$. This is quite an important case since the
fractional Sobolev spaces $W^{s,2}(\Omega)$ and $W^{s,2}_0(\Omega)$ turn out to be Hilbert spaces. They are usually denoted by $H^s(\Omega)$ and $H^s_0(\Omega)$, respectively.\\

$\bullet$ The spectral fractional Laplace operator: For $\Omega$ be a smooth bounded open subset of $\mathbb{R}^n$. The spectral fractional Laplace operator $\mathcal{A}^s$ is defined as follows. Let $\varphi_k$ be an eigenfunction of $-\Delta$ given by
\begin{equation}
\left\{\begin{array}{ccll}
-\Delta\varphi_k&=&\lambda_k\varphi_k & {\rm in} \ \ \Omega\; \\
\varphi_k&=& 0 & {\rm on} \ \ \partial\Omega
\end{array}\right.  ,
\end{equation}
where $\lambda_k$ is the corresponding eigenvalue of $\varphi_k,0<\lambda_1<\lambda_2\leq\lambda_3\leq\cdots\leq\lambda_k\rightarrow +\infty.$ Then, $\{\varphi_k\}_{k=1}^{\infty}$ is an orthonormal basic of $L^2(\Omega)$ satisfying 
\[
\int\limits_{\Omega}\varphi_j\varphi_kdx=\delta_{j,k}.
\]
We define the operator $\mathcal{A}^s$ for any $u\in C^\infty_0(\Omega)$ by

\begin{equation}\label{lapla}
\mathcal{A}^su=\sum_{k=1}^{\infty}\lambda_k^s\xi_k\varphi_k,
\end{equation}
where
\[
u=\sum_{k=1}^{\infty}\xi_k\varphi_k\text{  and  }\xi_k=\int\limits_{\Omega}u\varphi_kdx.
\]

$\bullet$ The restricted fractional Laplace operator: In this case we materialize the zero Dirichlet condition by restricting the operator to act only on functions that are zero outside $\Omega$. We will call the operator defined in such a way the restricted fractional Laplace operator. So defined, $(-\Delta)^s$ is a self-adjoint operator on $L^{2}(\Omega)$, with a discrete spectrum: we will denote by $\mu_{k}>0$, $k= 1, 2,\ldots$ its eigenvalues written in increasing order and repeated according to their multiplicity and we will denote by $\{\psi_{k}\}_k$ the corresponding set of eigenfunctions, normalized in $L^2(\Omega)$, where $\psi_{k}\in H^{2s}_0(\Omega)$. Eigenvalues $\mu_{k}$ (including multiplicities) satisfy
\[
0<\mu_{1}<\mu_{2}\leq\mu_{3}\leq\cdots\leq\mu_{k}\rightarrow +\infty.
\]

The spectral fractional Laplace operator $\mathcal{A}^s$ is related to (but different from) the restricted fractional Laplace operator $(-\Delta)^s$.

\begin{teor} The operators $(-\Delta)^s$ and $\mathcal{A}^s$ are not the same, since they have different eigenvalues and eigenfunctions. More precisely:
\begin{itemize}
    \item[(i)] the first eigenvalues of $(-\Delta)^s$ is strictly less than the one of $\mathcal{A}^s$.
    \item[(ii)] the eigenfunctions of $(-\Delta)^s$ are only H\"{o}lder continuous up to the boundary, differently from the ones of $\mathcal{A}^s$ that are as smooth up the boundary as the boundary allows.
\end{itemize}
\end{teor}
\n {\bf Proof.} See \cite{rafaella}.\; \fim\\

$\bullet$ Common notation. In the sequel we use $\mathcal{L}$ to refer to any of the two types of operators $\mathcal{A}^s$ or $(-\Delta)^s$, $0 < s < 1$. Each one is defined on a Hilbert space
\[
\Theta^{s}(\Omega)=\{u=\sum_{k=1}^{\infty}u_k\psi_{k}\in L^2(\Omega)\mid\sum_{k=1}^{\infty}\mu_{k}\vert u_k\vert^2<+\infty\}
\]
with values in its dual $\Theta^{s}(\Omega)'$. The Spectral Theorem allows to write $\mathcal{L}$ as
\[
\mathcal{L}u=\sum_{k=1}^{\infty}\mu_{k}u_k\psi_{k}
\]
for any $u\in \Theta^{s}(\Omega)$. Thus the inner product of $\Theta^{s}(\Omega)$ is given by
\[
\langle u,v\rangle_{\Theta^{s}(\Omega)}=\int\limits_{\Omega}\mathcal{L}^{1/2}u\mathcal{L}^{1/2}vdx=\int\limits_{\Omega}u\mathcal{L}vdx=\int\limits_{\Omega}v\mathcal{L}udx.
\]
We denote by $\Vert\cdot\Vert_{\Theta^{s}(\Omega)}$ the norm derived from this inner product. The notation in the formula copies the one just used for the second operator. When applied to the first one we put here $\psi_{k} = \varphi_k$, and $\mu_{k} = \lambda_{k}^s$. Note that $\Theta^{s}(\Omega)$ depends in principle on the type of operator and on the exponent $s$. It turns out that $\Theta^{s}(\Omega)$ independent of operator for each $s$, see \cite{sire}. We remark that $\Theta^{s}(\Omega)'$ can be described as the completion of the finite sums of the form
\[
f=\sum_{k=1}^{\infty}c_k\psi_{k} 
\]
with respect to the dual norm
\[
\Vert f\Vert_{\Theta^{s}(\Omega)'}=\sum_{k=1}^{\infty}\mu_{k}^{-1}\vert c_k\vert^2=\Vert\mathcal{L}^{-1/2}f\Vert_{L^2(\Omega)}^2=\int\limits_{\Omega}f\mathcal{L}^{-1}fdx
\]
and it is a space of distributions. Moreover, the operator $\mathcal{L}$ is an isomorphism between $\Theta^{s}(\Omega)$ and $\Theta^{s}(\Omega)'\simeq\Theta^{s}(\Omega)$, given by its action on the eigenfunctions. If $u,v\in\Theta^{s}(\Omega)$ and $f = \mathcal{L}u$ we have, after this isomorphism,
\[
\langle f, v\rangle_{\Theta^{s}(\Omega)'\times\Theta^{s}(\Omega)} = \langle u, v\rangle_{\Theta^{s}(\Omega)\times\Theta^{s}(\Omega)} =\sum_{k=1}^{\infty}\mu_{k}u_kv_k.
\]
If it also happens that $f\in L^{2}(\Omega)$, then clearly we get 
\[
\langle f, v\rangle_{\Theta^{s}(\Omega)'\times\Theta^{s}(\Omega)}=\int\limits_{\Omega}fv dx.
\]
We have $\mathcal{L}^{-1}:\Theta^{s}(\Omega)'\rightarrow\Theta^{s}(\Omega)$ can be written as
\[
\mathcal{L}^{-1}f(x)=\int\limits_{\Omega}G_{\Omega}(x,y)f(y)dy,
\]
where $G_{\Omega}$ is the Green function of operator $\mathcal{L}$ (see \cite{Green, ros176}). It is known that
\[
\Theta^s(\Omega)=\left\{
\begin{array}{llll}
L^2(\Omega) & {\rm if} \ \ s=0\\
H^{s}(\Omega)=H^{s}_0(\Omega) & {\rm if} \ \ s\in(0,\frac{1}{2})\\
H^{\frac{1}{2}}_{00}(\Omega) & {\rm if} \ \ s=\frac{1}{2}\\
H^{s}_0(\Omega) & {\rm if} \ \ s\in(\frac{1}{2},1]\\
H^{s}(\Omega)\cap H^{1}_0(\Omega) & {\rm if} \ \ s\in(1,2]
\end{array},
\right.
\]
where $H^{\frac{1}{2}}_{00}(\Omega):=\{u\in H^{1/2}(\Omega)\mid\int_{\Omega}\frac{u^2(x)}{d(x)}dx<+\infty\}.$\\

The next theorem gives a relation between the spectral fractional Laplace operator $\mathcal{A}^s$ and the restricted fractional Laplace operator $(-\Delta)^s$.

\begin{teor}
For $u\in H^{s}(\mathbb{R}^n)$, $u\geq 0$ and $\hbox{supp}(u)\subset\overline{\Omega}$, the following relation holds in the sense of distributions:
\[
\mathcal{A}^su\geq (-\Delta)^su.
\]
If $u\neq 0$ then this inequality holds with strict sign.
\end{teor}
\n {\bf Proof.} See \cite{compare}.\; \fim\\

In order to prove our main theorem, we should first introduce the concept of weak solution in both cases.

By weak solutions, we mean the following: Let $f\in L^{\frac{2n}{n+2s}}(\Omega)$. Given the problem
\begin{equation}\label{prob2}
\left\{\begin{array}{ccll}
\mathcal{A}^{s}u&=& f & {\rm in} \ \ \Omega\; \\
u&=& 0 & {\rm on} \ \ \partial\Omega
\end{array}\right.  
\end{equation}
we say that a function $u\in\Theta^s(\Omega)$ is a weak solution of (\ref{prob2}) provided
\[
\int\limits_{\Omega}\mathcal{A}^{s/2}u\mathcal{A}^{s/2}\phi dx=\int\limits_{\Omega}f\phi dx
\]
for all $\phi\in\Theta^s(\Omega)$. Now given the problem
\begin{equation}\label{prob1}
\left\{\begin{array}{ccll}
(-\Delta)^{s}u&=& f & {\rm in} \ \ \Omega\; \\
u&=& 0 & {\rm in} \ \ \mathbb{R}^n\setminus\Omega
\end{array}\right.  
\end{equation}
we say that a function $u$ is a weak solution of (\ref{prob1}) if $u\in H^s_{0}(\Omega)$, and
\[
\int\limits_{\mathbb{R}^n}(-\Delta)^{s/2}u(-\Delta)^{s/2}\varphi dx=\int\limits_{\Omega}f\varphi dx
\]
for all $\varphi\in H^{s}_{0}(\Omega)$. \\

\section{Solutions of system \ref{1}}

We will the proof for the spectral fractional Laplace operator $\mathcal{A}^s$. Similarly, follows the results for the restricted fractional Laplace operator $(-\Delta)^s$, changing the corresponding space. 

The existence result follows by applying the proof of [\cite{vander}, Theorem 2] for the case $s = 1$ with only minor
modifications.

We define the product Hilbert spaces
\[
E^{\alpha}(\Omega)=\Theta^{\alpha}(\Omega)\times\Theta^{2s-\alpha}(\Omega),\text{   }0<\alpha<2s
\]
where your inner product is given by
\[
\langle(u_1,v_1),(u_2,v_2)\rangle_{E^{\alpha}(\Omega)}=\langle\mathcal{A}^{\alpha/2}u_1,\mathcal{A}^{\alpha/2}u_2\rangle_{L^{2}(\Omega)}+\langle\mathcal{A}^{s-\alpha/2}v_1,\mathcal{A}^{s-\alpha/2}v_2\rangle_{L^{2}(\Omega)}.
\]
Let us remember that for the spectral fractional Laplace operator $\mathcal{A}^s$
\[
\Theta^{\alpha}(\Omega)=\{u=\sum_{k=1}^{\infty}u_k\varphi_{k}\in L^2(\Omega)\mid\sum_{k=1}^{\infty}\lambda_{k}^\alpha\vert u_k\vert^2<+\infty\}.
\]
We have $E^{\alpha}(\Omega)\hookrightarrow L^{q+1}(\Omega)\times L^{p+1}(\Omega)$ is compact for all exponents $q$ and $p$ satisfying
\[
q+1<\frac{2n}{n-2\alpha}\text{ and }p+1<\frac{2n}{n+2\alpha-4s}.
\]
We also have $\mathcal{A}^s:\Theta^\alpha(\Omega)\rightarrow\Theta^{\alpha-2s}(\Omega)$ is an isomorphism, see \cite{vander}. Hence
$$\left(\begin{array}{ccllrr}
0& \mathcal{A}^{s}  \\
\mathcal{A}^{s} & 0 
\end{array}\right): E^\alpha(\Omega)\rightarrow\Theta^{-\alpha}\times\Theta^{\alpha-2s}(\Omega)=E^\alpha(\Omega)'$$
is an isometry. We consider the Lagrangian
\begin{eqnarray*}
\mathcal{J}(u,v)&=&\int\limits_{\Omega}\mathcal{A}^{s/2}u\mathcal{A}^{s/2}vdx-\frac{\lambda}{2}\int\limits_{\Omega}u^2dx-\frac{\mu}{2}\int\limits_{\Omega}v^2dx\\
&-&\frac{1}{p+1}\int\limits_{\Omega}\vert v\vert^{p+1}dx-\frac{1}{q+1}\int\limits_{\Omega}\vert u\vert^{q+1}dx,
\end{eqnarray*}
i.e., a strongly indefinite functional. The o Hamiltonian is given by
\begin{equation}\label{2.2}
\mathcal{H}(u,v)=\frac{1}{p+1}\int\limits_{\Omega}\vert v\vert^{p+1}dx+\frac{1}{q+1}\int\limits_{\Omega}\vert u\vert^{q+1}dx.
\end{equation}
The quadratic part can again be written as
\[
A(u,v)=\frac{1}{2}\langle L(u,v),(u,v)\rangle_{E^{\alpha}(\Omega)}=\int\limits_{\Omega}\mathcal{A}^{s/2}u\mathcal{A}^{s/2}vdx-\frac{\lambda}{2}\int\limits_{\Omega}u^2dx-\frac{\mu}{2}\int\limits_{\Omega}v^2dx,
\]
where
$$L=\left(\begin{array}{ccllrr}
-\lambda\mathcal{A}^{-\alpha}& \mathcal{A}^{s-\alpha}  \\
\mathcal{A}^{\alpha-s} & -\mu\mathcal{A}^{\alpha-2s} 
\end{array}\right)$$
is bounded and self-adjoint. Unlike $L$ is not an isometry. 

In order to determine the spectrum of $L$, we note that $E^{\alpha}(\Omega)$ is the direct Hilbert space sum of the spaces $E_k, k=1,2,3,\ldots,$ where $E_k$ is the two-dimensional subspace of $E^{\alpha}(\Omega)$, spanned by $(\varphi_k,0)$ and $(0,\varphi_k)$. An orthonormal basis of $E_k$ is given by
\[
\left\{\frac{1}{\sqrt{2}}(\lambda_k^{-\alpha/2}\varphi_k,0),\frac{1}{\sqrt{2}}(0,\lambda_k^{\alpha/2-s}\varphi_k)\right\}.
\]
Every $E_k$ is invariant under $L$, and in $E_k$ the restriction of $L$ is given the symmetric matrix
$$L^k=\left(\begin{array}{ccllrr}
-\lambda\lambda_k^{-\alpha}& 1  \\
1 & -\mu\lambda_k^{\alpha-2s} 
\end{array}\right).$$
The eigenvalues of $L^k$ are given by
\[
\nu_k^\pm=-\frac{\lambda\lambda_k^{-\alpha}+\mu\lambda_k^{\alpha-2s}}{2}\pm\sqrt{\left(\frac{\lambda\lambda_k^{-\alpha}+\mu\lambda_k^{\alpha-2s}}{2}\right)^{2}+1-\frac{\lambda\mu}{\lambda_k^{2s}}},
\]
with corresponding eigenvectors
\[
\left(1,\frac{\lambda\lambda_k^{-\alpha}-\mu\lambda_k^{\alpha-2s}}{2}\pm\sqrt{\left(\frac{\lambda\lambda_k^{-\alpha}-\mu\lambda_k^{\alpha-2s}}{2}\right)^{2}+1}\right).
\]
We have $\nu_{k}^{-}<0<\nu_{k}^+$ if $\lambda\mu<\lambda_{k}^{2s}$. If $\lambda\mu>\lambda_{k}^{2s}$ the signs of $\nu_{k}^+$ and $\nu_{k}^-$ are the same: positive (negative) if $\lambda$ and $\mu$ are negative (positive). If $\lambda\mu=\lambda_{k}^{2s}$, then $\nu_k^+=0\text{ }(\nu_{k}^-=0)$ if $\lambda$ and $\mu$ are positive (negative). Also note that
\begin{equation}\label{4.7}
\nu_k^\pm\rightarrow\pm1\text{ as }k\rightarrow\infty.
\end{equation}
Let $E^+\text{ }(E^-)$ be the subspace spanned by eigenvectors with positive (negative) eigenvalues, and $E^0$ the nullspace of $L$. Then
\begin{equation}\label{4.8}
E^{\alpha}(\Omega)=E^+\oplus E^-\oplus E^0.
\end{equation}
It follows that both $E^+$ and $E^-$ are infinite dimensional, and that $E^0$ has finite dimension: $\lambda\mu\neq\lambda_{k}^{2s}$ implies $\dim E^0=0$ while for $\lambda\mu=\lambda_{k}^{2s}$ the dimension of $E^0$ is equal to the multiplicity of $\lambda_k^s$.
We introduce a equivalent (inner product) norm $\Vert\cdot\Vert_\ast$ on $E^{\alpha}(\Omega)$ by
\begin{equation}\label{4.9}
\langle L(u,v)^+,(u,v)^+\rangle_{E^{\alpha}(\Omega)}-\langle L(u,v)^-,(u,v)^-\rangle_{E^{\alpha}(\Omega)}+\Vert (u,v)^0\Vert_{L^2(\Omega)\times L^2(\Omega)}^2=\Vert(u,v)\Vert_\ast^2.
\end{equation}
Note that the equivalence of (\ref{4.9}) and $\Vert\cdot\Vert_{E^\alpha(\Omega)}$ follows from (\ref{4.7}) and the fact $E^0$ is finite dimensional.

The derivative of $A(u,v)$ defines a bilinear form
\begin{equation}\label{1.27}
B((u_1,v_1),(u_2,v_2))=A'(u_1,v_1)(u_2,v_2)=\langle L(u_1,v_1),(u_2,v_2)\rangle_{E^{\alpha}(\Omega)},
\end{equation}
where $(u_1,v_1),(u_2,v_2)\in E^{\alpha}(\Omega)$ with
\begin{equation}\label{1.28}
A(u_1,v_1)=\frac{1}{2}B((u_1,v_1),(u_1,v_1))\text{ and }B((u_1,v_1)^+,(u_1,v_1)^-)=0.
\end{equation}

The proof of Theorem \ref{existhamiltinian} is based on an application of the following result of Benci and Rabinowitz \cite{benci}.

\begin{teor} \label{benci}(Indefinite Functional Theorem). Let $H$ be a real Hilbert sapce with $H = H_1\oplus H_2$. Suppose $\mathcal{J}\in C^{1}(H,\mathbb{R})$ satisfies the Palais-smale condition, and
    \begin{itemize}
      \item[(i)] $\mathcal{J}(w) =\frac{1}{2}(Lw, w)_H-\mathcal{H}(w)$, where $\mathcal{J}: H \rightarrow H$ is bounded and self-adjoint, and $L$ leaves $H_1$ and $H_2$ invariant.
      \item[(ii)] $\mathcal{H}'$ is compact.
      \item[(iii)] there exists a subspace $\overline{H}\subset H$ and sets $S\subset H, Q\subset\overline{H}$ and constants $\beta > \omega$ such that
\begin{itemize}
\item[(1)] $S\subset H_1$ and $\mathcal{J}\mid_{S}\geq\beta$.
\item[(2)] $Q$ is bounded and $\mathcal{J}\leq\omega$ on the boundary $\partial Q$ of $Q$ in $\overline{H}$.
\item[(3)] $S$ and $\partial Q$ link.
\end{itemize}
     \end{itemize}
Then $\mathcal{J}$ possesses a critical value $c\geq\beta$.
\end{teor}

{\bf Proof of Theorem \ref{existhamiltinian}.} We apply Theorem \ref{benci} with the spaces $H = E^{\alpha}(\Omega)$, $H_1 = E^+$, and $H_2 =E^0\oplus E^-$. Apart from the choice of $\alpha$, this is standard, and follows from condition \ref{5}. We will use $\textbf{u}=(u,v)\in E^\alpha(\Omega)$ to simplify the notation of this proof.

To show that $\mathcal{J}$ is continuously differentiable and satisfies $(i)$ and $(ii)$, it suffices to observe that $\mathcal{H}(\textbf{u})$, defined by (\ref{2.2}), is continuously differentiable with compact derivative $\mathcal{H}'(\textbf{u})$. 

If $\lambda\mu\neq\lambda_{k}^{2s}$ for all $k$ we have $\dim E^0=0$, and (\ref{4.9}) reduces to
\begin{equation}
\langle L(\textbf{u}^+),\textbf{u}^+\rangle_{E^{\alpha}(\Omega)}-\langle L(\textbf{u}^-),\textbf{u}^-\rangle_{E^{\alpha}(\Omega)}=\Vert\textbf{u}\Vert_\ast^2.
\end{equation}

Step 1. The geometry conditions $(iii)(1),(2),(3)$ are satisfied. Let $\rho$, $t_1>\rho$ and $t_2$ be positive numbers to be specified later on, and we take for $\textbf{e}^+$ an eigenvector in $E^+$, such that $\textbf{e}^+$ belongs to some $E_k$ with the other eigenvector $\textbf{e}^-$ in $E_k$ belonging to $E^-$ ($\textbf{e}^+$ and $\textbf{e}^-$ normalized with respect to $\Vert\cdot\Vert_\ast$). Note that $\textbf{e}^-$ is the only eigenvector of $L$ not perpendicular to $\textbf{e}^+$ in $L^2(\Omega)\times L^2(\Omega)$. We set $[0,t_1\textbf{e}^+]=\{t\textbf{e}^+\mid 0\leq t\leq t_1\}$, and
\[
Q=[0,t_1\textbf{e}^+]\oplus(\overline{B}_{t_2}\cap E^-),\text{  }\overline{H}=span[\textbf{e}^+]\oplus E^-,\text{  }S=\partial B_{\rho}\cap E^+,
\] 
where $B_{R}$ denotes an open ball with radius $R$ centered at the origin.

On $E^+$ the quadratic part $A(\textbf{u})$ of $\mathcal{J}$ reduces to $\frac{1}{2}\Vert\textbf{u}\Vert_\ast^2$, so that $A$ has a strict local minimum on $E^+$ at $\textbf{u}^+=0$. Thus the same is true for $\mathcal{J}$. Indeed there is $C>0$ such that
\[
\mathcal{J}(\textbf{u}^+)\geq\frac{1}{2}\Vert\textbf{u}^+\Vert_\ast^2-C\Vert\textbf{u}\Vert_\ast^{p+1}-C\Vert\textbf{u}\Vert_\ast^{q+1}
\]
for all $\textbf{u}^+\in E^+$. Thus we can fix $\rho>0$ and $\beta>0$ such that $\mathcal{J}(\textbf{u})\geq\beta$ on $S$, and $(iii)(1)$ is satisfied.

Next we show that for suitable choices of $t_1$ and $t_2$ the function $\mathcal{J}(\textbf{u})$ is nonpositive on $\partial Q$. Thus we prove that $(iii)(2)$ holds with $\omega=0$. Note that the boundary $\partial Q$ of the cylinder $Q$ is taken in the space $\overline{H}$, and consists of three parts, namely the bottom $Q\cap\{t=0\}$, the lid $Q\cap\{t=t_1\}$, and the lateral boundary $[0,t_1\textbf{e}^+]\oplus(\partial B_{t_2}\cap E^-)$. Clearly $\mathcal{J}(\textbf{u})\leq 0$ on the bottom because $A(\textbf{u})\leq 0$ in $E^-$, and the functional $\mathcal{H}(\textbf{u})$ is nonnegative. For the remaining two parts of the boundary we first observe that, for $\textbf{u}=\textbf{u}^-+t\textbf{e}^+\in\overline{H}$,
\[
\mathcal{J}(\textbf{u}^-+t\textbf{e}^+)=\frac{1}{2}t^2 -\frac{1}{2}\Vert\textbf{u}^-\Vert_\ast^2-\mathcal{H}(\textbf{u}^-+t\textbf{e}^+).
\]  
Since $p,q>1$ there is $\gamma>2$ such that $\gamma<p+1$ and $\gamma<q+1$. Thus there is $B_1>0$ such that 
\begin{equation}\label{2.7}
\mathcal{H}(\textbf{u}^-+t\textbf{e}^+)\geq B_1\int\limits_{\Omega}\vert\textbf{u}^-+t\textbf{e}^+\vert^\gamma.
\end{equation}
So, writing $\textbf{u}^-=h\textbf{e}^-+\textbf{u}_2^-$, where $h$ is a real number, and $\textbf{u}_2^-\in E^-$ is perpendicular to $\textbf{e}^-$ in $E^\alpha(\Omega)$, $\textbf{u}_2^-$ is perpendicular to $\textbf{e}^-$ and $\textbf{e}^+$ in $L^2\times L^2$, and we proceed from (\ref{2.7}) to conclude
\begin{eqnarray*}
\mathcal{H}(\textbf{u}^-+t\textbf{e}^+)&\geq & B_1\int\limits_{\Omega}\vert h\textbf{e}^-+\textbf{u}_2^-+t\textbf{e}^+\vert^\gamma\geq B_2\left(\int\limits_{\Omega}\vert h\textbf{e}^-+\textbf{u}_2^-+t\textbf{e}^+\vert^2\right)^{\gamma/2}\\
&=&B_2\left(\int\limits_{\Omega}\vert\textbf{u}_2^-\vert^2+\int\limits_{\Omega}\vert h\textbf{e}^-+t\textbf{e}^+\vert^2\right)^{\gamma/2}\geq B_2\left(\int\limits_{\Omega}\vert h\textbf{e}^-+t\textbf{e}^+\vert^2\right)^{\gamma/2}\\
&\geq &B_2\left(t^2\sin^2\theta\int\limits_{\Omega}\vert \textbf{e}^+\vert^2\right)^{\gamma/2}. 
\end{eqnarray*}
Here $\theta$ is the (positive) angle between $\textbf{e}^+$ and $\textbf{e}^-$ with respect to the inner product in $L^2(\Omega)\times L^2(\Omega)$. Thus
\begin{equation}\label{2.8}
\mathcal{H}(\textbf{u}^-+t\textbf{e}^+)\geq B_3t^\gamma.
\end{equation}
For $\mathcal{J}$ this yields
\begin{equation}\label{2.9}
\mathcal{J}(\textbf{u}^-+t\textbf{e}^+)\leq\frac{1}{2}t^2-B_3t^\gamma-\frac{1}{2}\Vert\textbf{u}^-\Vert_\ast^2.
\end{equation}
Now choose $t_1$, $t_2$ such that
\[
\psi(t)=\frac{1}{2}t^2-B_3t^\gamma\leq 0\text{  }\forall t\leq t_1,t_2^2>2\max_{t\geq 0}\psi(t)
\]
to make $\mathcal{J}$ negative on the lid and on the lateral boundary, respectively. Then $(iii)(2)$ is satisfied with $\omega=0$.

For the proof of $(iii)(3)$ we refer to \cite{benci}.\\

Step 2. The Palais-Smale condition. We have to show that any sequence $\{\textbf{u}_n\}$ in $E^\alpha(\Omega)$ satisfying 
\begin{equation}\label{2.11}
\mathcal{J}(\textbf{u}_n)\text{ is bounded in }E^\alpha(\Omega),\text{   }\mathcal{J}'(\textbf{u}_n)\rightarrow 0,\text{ as }n\rightarrow\infty
\end{equation} 
has a convergence subsequence. The key point here is to prove that such a sequence is necessarily bounded in $E^\alpha(\Omega)$. For then the compactness of $\mathcal{H}'$ implies, since $\mathcal{J}'(\textbf{u}_n)$ converges in $E^\alpha(\Omega)'$, that a subsequence of $A'(\textbf{u}_n)$ also converges. In view of (\ref{1.27}), we then also have $L\textbf{u}_n$ and $\textbf{u}_n$ converging in $E^\alpha(\Omega)$, because $L$ is invertible.

To prove that $\textbf{u}_n$ is bounded we proceed as follows. By (\ref{1.27}), (\ref{1.28}) and (\ref{2.11}), for some $M>0$, and arbitrarily small $\varepsilon>0$, omitting the subscripts, 
\begin{eqnarray*}
M+ \varepsilon\Vert\textbf{u}\Vert_\ast &\geq & \mathcal{J}(\textbf{u})-\frac{1}{2}\langle\mathcal{J}'(\textbf{u}),\textbf{u}\rangle\\
&=&\int\limits_{\Omega}\left(\frac{1}{2}u\vert u\vert^{q}-\frac{1}{q+1}\vert u\vert^{q+1}+\frac{1}{2}v\vert v\vert^{p}-\frac{1}{p+1}\vert v\vert^{p+1}\right)dx.
\end{eqnarray*}
Thus
\begin{equation}\label{2.12}
M+ \varepsilon\Vert\textbf{u}\Vert_\ast \geq\left(\frac{1}{2}-\frac{1}{\gamma}\right)\left[\int\limits_{\Omega}\vert u\vert^{q+1}dx+\int\limits_{\Omega}\vert v\vert^{p+1}dx\right] - C,
\end{equation}
where $C$ is a constant taking care of the $u,v$-values between $-\eta$ and $\eta$. Hence, for some new constants $C,\varepsilon>0$,
\begin{equation}\label{2.13}
\int\limits_{\Omega}\vert u\vert^{q+1}dx+\int\limits_{\Omega}\vert v\vert^{p+1}dx\leq C+\varepsilon\Vert\textbf{u}\Vert_\ast.
\end{equation}
Still omitting subscripts, and writing $\textbf{u}^\pm=(u^\pm,v^\pm)$, we also have, by (\ref{2.11}),
\begin{eqnarray*}
\Vert\textbf{u}^\pm\Vert_\ast^2 -\varepsilon\Vert\textbf{u}^\pm\Vert_\ast &\leq & \vert\langle L(\textbf{u}),\textbf{u}^\pm\rangle_{E^{\alpha}(\Omega)}-\langle \mathcal{J}'(\textbf{u}),\textbf{u}^\pm\rangle\vert=\vert\langle\mathcal{H}'(\textbf{u}),\textbf{u}^\pm\rangle\vert\\
&=&\left\vert\int\limits_{\Omega}\vert u\vert^{q}u^\pm dx +\int\limits_{\Omega} \vert v\vert^{p}v^\pm dx\right\vert\\
&\leq &\left(\int\limits_{\Omega}\vert u\vert^{q+1}dx\right)^{q/q+1}\Vert u^\pm\Vert_{L^{q+1}(\Omega)}+\left(\int\limits_{\Omega}\vert v\vert^{p+1}dx\right)^{p/p+1}\Vert v^\pm\Vert_{L^{p+1}(\Omega)}\\
&\leq & C\left[\left(\int\limits_{\Omega}\vert u\vert^{q+1}dx\right)^{q/q+1}+\left(\int\limits_{\Omega}\vert v\vert^{p+1}dx\right)^{p/p+1}\right]\Vert\textbf{u}^\pm\Vert_\ast.
\end{eqnarray*}
Dividing the first and the last expression by $\Vert\textbf{u}^\pm\Vert_\ast$ we obtain
\begin{equation}\label{2.14}
\Vert\textbf{u}^\pm\Vert_\ast-\varepsilon\leq C\left[\left(\int\limits_{\Omega}\vert u\vert^{q+1}dx\right)^{q/q+1}+\left(\int\limits_{\Omega}\vert v\vert^{p+1}dx\right)^{p/p+1}\right].
\end{equation}
Combining (\ref{2.14}) for $\textbf{u}=\textbf{u}^++\textbf{u}^-$, together with (\ref{2.13}), it follows that, possibly for some new constant,
\begin{equation}\label{2.15}
\Vert\textbf{u}\Vert_\ast\leq C\left[\left(C+\varepsilon\Vert\textbf{u}\Vert_\ast\right)^{q/q+1}+\left(C+\varepsilon\Vert\textbf{u}\Vert_\ast\right)^{p/p+1}\right],
\end{equation}
which keeps $\Vert\textbf{u}\Vert_\ast$ away from infinity. This implies that the Palais-Smale condition is satisfied.

If $\lambda\mu=\lambda_{k}^{2s}$ for some $k$, the proof is slightly more complicated. We define 
\[
Q=[0,t_1\textbf{e}^+]\oplus(\overline{B}_{t_2}\cap H_2),\text{  }\overline{H}=span[\textbf{e}^+]\oplus H_2,\text{  }S=\partial B_{\rho}\cap H_1,
\] 
Elements of $H_1=E^+$ are denoted by $\textbf{u}^+$, and elements of $H_2=E^0\oplus E^-$ by $\textbf{u}^0+\textbf{u}^-$.

To verify the geometric conditions in Step 1, we have to estimate
\[
\mathcal{J}(\textbf{u}^0+\textbf{u}^-+t\textbf{e}^+)=\frac{1}{2}t^2 -\frac{1}{2}\Vert\textbf{u}^-\Vert_\ast^2-\mathcal{H}(\textbf{u}^0+\textbf{u}^-+t\textbf{e}^+)
\]  
on the boundary of the cylinder $Q$. The lateral boundary however is no longer given by $\Vert\textbf{u}^-\Vert_\ast=t_2$ but by $\Vert\textbf{u}_0+\textbf{u}^-\Vert_\ast=t_2$. Thus if $t_2$ is large, the norm $\Vert\textbf{u}^-\Vert_\ast$ can still be small on the lateral boundary, provided $\Vert\textbf{u}^0\Vert_\ast$ is large. Therefore we now estimate $\mathcal{H}(\textbf{u}^0+\textbf{u}^-+t\textbf{e}^+)$ from below by changing (\ref{2.8}) by
\begin{eqnarray*}
\mathcal{H}(\textbf{u}^0+\textbf{u}^-+t\textbf{e}^+)&\geq & B_2\left(\int\limits_{\Omega}\vert \textbf{u}^0+\textbf{u}_2^-+t\textbf{e}^+\vert^2\right)^{\gamma/2}\\
&\geq &B_2\left(\sin^2\theta\int\limits_{\Omega}\vert \textbf{u}^0+t\textbf{e}^+\vert^2\right)^{\gamma/2}\\
&=&B_2\sin^\gamma\theta(\Vert\textbf{u}^0\Vert_\ast^2+t^2)^{\gamma/2}.
\end{eqnarray*}
Here $\theta$ is the (positive) angle between $E^-$ and $E^0\oplus[\textbf{e}^+]$ with respect to the inner product in $L^2(\Omega)\times L^2(\Omega)$. Thus
\begin{equation}\label{4.13}
\mathcal{H}(\textbf{u}^0+\textbf{u}^-+t\textbf{e}^+)\geq B_4\Vert\textbf{u}^0\Vert_\ast^2+B_5t^\gamma.
\end{equation}
The analog of (\ref{2.9}) is
\begin{equation}\label{4.14}
\mathcal{J}(\textbf{u}^0+\textbf{u}^-+t\textbf{e}^+)\leq\frac{1}{2}t^2-B_5t^\gamma-B_6(\Vert\textbf{u}^-\Vert_\ast^2+\Vert\textbf{u}^0\Vert_\ast^2)=\frac{1}{2}t^2-B_5t^\gamma-B_6(\Vert\textbf{u}^-+\textbf{u}^0\Vert_\ast^2).
\end{equation}
The proof now proceeds along the same lines as before.

It remains to show that $\mathcal{J}$ satisfies the Palais-Smale condition. So let $\{\textbf{u}_n\}$ in $E^\alpha(\Omega)$ be a sequence with $\mathcal{J}(\textbf{u}_n)$ bounded in $E^\alpha(\Omega)$, and $\mathcal{J}'(\textbf{u}_n)\rightarrow 0$. We have to do some extra work to show that such a sequence is bounded. Using the decomposition (\ref{4.8}) estimates (\ref{2.13}) and (\ref{2.14}) for $\textbf{u}^\pm$ remain the same. Thus we have, instead of (\ref{2.15}),
\begin{equation}\label{4.15}
\Vert\textbf{u}^\pm\Vert_\ast\leq C\left[\left(C+\varepsilon\Vert\textbf{u}\Vert_\ast\right)^{q/q+1}+\left(C+\varepsilon\Vert\textbf{u}\Vert_\ast\right)^{p/p+1}\right].
\end{equation}
To control the component $\textbf{u}^0$ we modify (\ref{2.12}) and derive
\begin{eqnarray*}
M+ \varepsilon\Vert\textbf{u}\Vert_\ast &\geq & \int\limits_{\Omega}\left(\frac{1}{2}u\vert u\vert^{q}-\frac{1}{q+1}\vert u\vert^{q+1}+\frac{1}{2}v\vert v\vert^{p}-\frac{1}{p+1}\vert v\vert^{p+1}\right)dx\\
&\geq &\left(\frac{\gamma}{2}-1\right)\mathcal{H}(\textbf{u})\geq B_1\int\limits_{\Omega}\vert \textbf{u}\vert^{\gamma}\geq B_2\left(\int\limits_{\Omega}\vert \textbf{u}^-+\textbf{u}^0+\textbf{u}^+\vert^{2}\right)^{\gamma/2}\\
&\geq &B_2\left(\sin^2\theta\int\limits_{\Omega}\vert \textbf{u}^0\vert^2\right)^{\gamma/2},
\end{eqnarray*}
where $\theta$ is the (positive) angle between $E^0$ and $E^+\oplus E^-$ with respect to the inner product in $L^2(\Omega)\times L^2(\Omega)$. Thus
\begin{equation}\label{4.16}
M+ \varepsilon\Vert\textbf{u}\Vert_\ast \geq B_3\Vert\textbf{u}^0\Vert_\ast^\gamma.
\end{equation}
Combining (\ref{4.15}) and (\ref{4.16}) we obtain
\begin{equation}
\Vert\textbf{u}\Vert_\ast\leq C\left[\left(C+\varepsilon\Vert\textbf{u}\Vert_\ast\right)^{q/q+1}+\left(C+\varepsilon\Vert\textbf{u}\Vert_\ast\right)^{p/p+1}+\left(C+\varepsilon\Vert\textbf{u}\Vert_\ast\right)^{1/\gamma}\right],
\end{equation}
and as before this implies that the sequence $\textbf{u}_n$ is bounded. This concludes the part of existence of the theorem.

Now the part of regularity of the theorem. First we shall prove the $L^{\infty}$ estimate of Brezis-Kato type.

\begin{propo} Suppose that $\lambda$ and $\mu$ are fixed real numbers and $p,q>1$ satisfies (\ref{5}). Let $(u, v)$ be a weak solution of (\ref{1}). Then we have $u\in L^{\infty}(\Omega)$ and $v\in L^\infty(\Omega)$.
\end{propo}

{\bf Proof.} Letting $a(x)=(\mu+ v^{p-1}(x))$ and $b(x)=(\lambda+ u^{q-1}(x))$, we have $a\in L^{\frac{p+1}{p-1}}(\Omega)$ and $b\in L^{\frac{q+1}{q-1}}(\Omega)$. Now we write (\ref{1}) as
\begin{equation}
\left\{
\begin{array}{llll}
\mathcal{L}u = a(x)v & {\rm in} \ \ \Omega\\
\mathcal{L}v = b(x)u & {\rm in} \ \ \Omega\\
\end{array}.
\right.
\end{equation}

Then one can follow the proof of [\cite{choi}, Proposition 5.2] if $\mathcal{L}=\mathcal{A}^s$ and the proof of [\cite{EM1}, Proposition 3.1] if $\mathcal{L}=(-\Delta)^s$ for conclude the result. \; \fim\\

From the above Proposition, we have $\lambda u+u^q,\mu v+v^p\in L^{\infty}(\Omega)$. If $\mathcal{L}=\mathcal{A}^s$ by regularity result (see \cite{cabre1} or Proposition 3.1 of \cite{T}) we have $v,u\in C^{2s}(\overline{\Omega})$. Hence it holds that $\lambda u+u^q,\mu v+v^p \in C^{2s}(\overline{\Omega})$. Again, we can apply regularity result to deduce that $u\in C^{4s}(\overline{\Omega})$. Iteratively, we can raise the regularity so that $(u,v) \in C^{1,\eta}(\overline{\Omega}) \times C^{1,\eta}(\overline{\Omega})$ for some $\eta \in (0, 1)$. Thus by Theorem 1.1 of \cite{yang} it follows that $(u,v)\in (C^{2,\eta}(\Omega)\cap C^{1,\eta}(\overline{\Omega}))^2$. 

Finally for the case $\mathcal{L}=(-\Delta)^s$, $C^\eta$ regularity of $u$ and $v$ in $\R^n$ for some $\eta \in(0,1)$ is obtained from each equation by evoking Proposition 1.4 of \cite{ROS}.

\end{document}